\documentclass{gtart}

\usepackage{amsmath,amsfonts,amscd,flafter,rlepsf}

\newtheorem{theorem}{Theorem}[section]

\newtheorem{lemma}[theorem]{Lemma}

\theoremstyle{remark}

\newtheorem{remark}[theorem]{Remark}

\newcommand\Rk{\mathrm{rk}}

\newcommand{\HF}{HF}

\newcommand{\Z}{\mathbb{Z}}

\newcommand{\cm}{\cdot}

\newcommand{\ModSWfour}{\mathcal{M}}
\newcommand{\ModFlow}{\ModSWfour}

\newcommand{\SpinC}{{\mathrm{Spin}}^c}

\newcommand\Hom{\mathrm{Hom}}

\newcommand\abuts\Rightarrow
\newcommand\Sym{\mathrm{Sym}}

\newcommand\ModSphere{\ModFlow\left({\mathbb S}\longrightarrow 
\Sym^{g-1}(\Sigma_{1})\times \Sym^2(\Sigma_{2})\right)}
\newcommand\ModSpheres\ModSphere

\newcommand\HFred{\HF_{\rm red}}

\newcommand\HFp{\HFb}

\newcommand\HFb{HF^+}

\newcommand\UnparModSp{\widehat \ModSp}
\newcommand\UnparModFlow\UnparModSp
\newcommand\Mod\ModSp

\newcommand\PD{\mathrm{PD}}

\newcommand{\spinct}{\mathfrak t}

\newcommand\ModMaps{\mathcal M}
\newcommand\ModSp\ModMaps

\newcommand\Dual{\mathcal D}
\newcommand\Duality\Dual

\newcommand\Char{\mathrm{Char}}
\newcommand\Combp{{\mathbb H}^+}
\newcommand\DCombp{{\mathbb K}^+}

\newcommand\InjMod[1]{{\mathcal T}^+_{#1}}

\begin{document}

\title{Calculation of Heegaard Floer homology for a class of Brieskorn spheres}

\author{Raif Rustamov}
\address{The Program in Applied and Computational Mathematics, Princeton University\\New Jersey 08540, USA}
\email{rustamov@princeton.edu}

\begin{abstract}
In this paper we calculate Ozsv\'ath-Szab\'o Floer homology group $\HFp$ for Brieskorn spheres $\Sigma(2,2n+1,4n+3)$.
\end{abstract}
\maketitle

\section {Introduction}
In \cite{HolDisk} and \cite{HolDiskFour} topological invariants for closed oriented three manifolds and cobordisms between them are defined by using a construction from symplectic geometry. The resulting Floer homology package has many properties of a topological quantum field theory.

The construction of Heegard Floer homology is more combinatorial in flavor than the corresponding gauge theoretical constructions of Donaldson-Floer (see~\cite{DonaldsonFloer}) and Seiberg-Witten theories (see~\cite{Witten}, \cite{Morgan}, \cite{KMthom}).
However, the construction still depends on profoundly analytic objects - holomorphic disks. As a result, one is tempted to consider classes of manifolds which allow a completely combinatorial description of their Heegaard Floer homologies.

In \cite{Plumbing} a class of plumbing three-manifolds is studied. It is proved there that $\HFp$ of these manifolds can be expressed in terms of equivariant maps. Based on this we calculate Heegaard Floer homology of a class of Brieskorn spheres.  

Let $\InjMod{0}$ denote the $\Z[U]$-module which is the quotient of
$\Z[U,U^{-1}]$ by the submodule $U\cm \Z[U]$, graded so that the element $U^{-d}$ (for $d\geq 0$) is supported in degree $2d$. For a positive integer $r$ define $\Z^r_{(k)}$ to be the quotient of $\Z[U]$ by $U^r\Z[U]$, where $U^{r-1}$ lies in degree $k$, and multiplication by $U$ decreases grading by $2$.

Let $q_i = i(i+1)$, and $p_i$ with $i=1,2,...$ denote the $i$'th element of the sequence $$1,1,2,2,3,3,4,4,....$$

There is a unique $\SpinC$ structure on $Y=-\Sigma(2, 2n+1,4n+3)$ and we will supress it from the notation.

Our main result is the following identification of $\HFp(-\Sigma(2, 2n+1, 4n+3))$.
\begin{theorem}
\label{theorem:main}
For any positive integer $n$ we have
$$\HFp(-\Sigma(2, 2n+1,4n+3)) \cong \InjMod{0} \oplus \Z^{p_n}_{(0)}\oplus  \bigoplus_{i=1}^{n-1} (\Z^{p_{i}}_{(q_{n-i})} \oplus \Z^{p_{i}}_{(q_{n-i})})$$
\end{theorem}

It would be interesting to compare this calculation with the corresponding
analogues in instanton Floer homology (see for
example~\cite{FintushelStern}) and Seiberg-Witten
theory (see for example~\cite{MOY}).
\begin{remark}
In this paper we study $\HFp$ for a special infinite family of Brieskorn spheres. For more computations see also \cite{Nemethi}.
\end{remark}
\medskip
\noindent{\bf{Acknowledgments}}\qua No words can express my gratitude to my advisor  Zolt\'an Szab\'o who kindly and patiently taught almost everything I know on Heegaard Floer homology and beyond. I would like to thank Paul Seymour for his invaluable support and encouragement. Many thanks go to Eaman Eftekhary and Andr\'as N\'emethi for helpful conversations.

\section{Heegaard Floer homology of a plumbing}
According to the combinatorial description of Heegard Floer homology groups (for a quite large class of plumbed three-manifolds) given in \cite{Plumbing} there are two main steps in the calculation: finding the basic vectors and finding the minimal relationships between them. The following is a review of the background information and algorithms needed.  

Let G be a weighted graph and let $m(v)$ and $d(v)$ be respectively the weight and the degree of the vertex $v$.
We denote by $X(G)$ the four-manifold with boundary having $G$ as its plumbing diagram. Let $Y(G)$ be the oriented three-manifold which is the boundary of $X(G)$
.

For $X=X(G)$, the group $H_2(X;\Z)$ is the lattice freely spanned by vertices of $G$.
 Denoting by $[v]$ the homology class in $H_2(X;\Z)$ corresponding to the vertex $v$ of $G$, the values of the intersection form of $X$ on the basis are given by $[v]\cm[v] = m(v)$; 
 $[v]\cm[w] = 1$ if $vw$ is an edge of $G$ and $[v]\cm[w]=0$ otherwise. The graph $G$ is called \emph{negative-definite} if the form is negative-definite. A vertex $v$ is said to be a \emph{bad vertex} of $G$ if $m(v)>-d(v)$.

Denoting by $\Char(G)$  the set of characteristic vectors for the
intersection form define 
$$\Combp(G)\subset \Hom(\Char(G),\InjMod{0})$$
to be the set of finitely supported functions satisfying the following relations for all characteristic vectors $K$ and vertices $v$:

\begin{equation}
\label{eq:AdjRel}
U^n\cm \phi(K+2\PD[v]) = \phi(K),
\end{equation}
if $2n=\langle K,v \rangle + v\cm v \geq 0$ and 
\begin{equation}
\label{eq:AdjRel2}
\phi(K+2\PD[v]) = U^{-n}\cm \phi(K)
\end{equation}
for $n<0$.

A grading on $\Combp(G)$ is introduced as follows: we say that an element $\phi\in \Combp(G)$ is homogeneous of degree $d$ if 
for each characteristic vector $K$ with $\phi(K)\neq 0$, $\phi(K)\in \InjMod{0}$ is a homogeneous element with:
\begin{equation}
\label{eq:DefOfDegree}
\deg(\phi(K))-\left(\frac{K^2+|G|}{4}\right)=d.
\end{equation}

After decomposing $\Combp(G)$ according to the $\SpinC$ structures over $Y$he following theorem is proved in \cite{Plumbing}:

\begin{theorem}
\label{intro:SomePlumbings}
Let $G$ be a negative-definite weighted forest with at most one bad vertex. Then, for each $\SpinC$ structure
$\spinct$ over $-Y(G)$, there is an isomorphism of graded
$\Z[U]$ modules,
$$\HFp(-Y(G),\spinct)\cong \Combp(G,\spinct).$$
\end{theorem}

For calculational purposes it is helpful to adopt a dual point of view. Let $\DCombp(G)$ be the set of equivalence
 classes of elements of $\Z^{\geq 0} \times \Combp(G)$ (and we write $U^m\otimes K$ for the pair $(m, K)$) 
under the following equivalence relation.
For any vertex $v$ let
$$2n=\langle K,v \rangle + v\cm v.$$
If $n\geq 0$, then
\begin{equation}
\label{equation:rel1}
U^{n+m}\otimes (K+2\PD[v]) \sim U^m\otimes K,
\end{equation}
while if $n\leq 0$, then
\begin{equation}
\label{equation:rel2}
U^m\otimes (K+2\PD[v]) \sim U^{m-n}\otimes K.
\end{equation}

Starting with a map $$\phi\colon \Char(G)\longrightarrow \InjMod{0},$$
consider an induced map $${\widetilde \phi}\colon \Z^{\geq 0}\times
\Char(G)\longrightarrow \InjMod{0}$$
defined by $${\widetilde \phi}(U^n\otimes K)=U^n\cm \phi(K).$$
Clearly, the set of finitely-supported functions $\phi\colon \Char(G)\longrightarrow
\InjMod{0}$ whose induced map ${\widetilde \phi}$ descends to
$\DCombp(G)$ is precisely $\Combp(G)$.

A \emph{basic element} of $\DCombp(G)$ is one whose equivalence class does not contain  any of $U^m\otimes K$ with $m>0$.
Given two non-equivalent basic elements $K_1 = U^0 \otimes K_1$ and $K_2 = U^0 \otimes K_2$ in the same $\SpinC$ structure, one can find positive
 integers $n$ and $m$ such that 
$$ U^{n} \otimes K_1 \sim U^{m} \otimes K_2.$$ 
If, moreover, the numbers $n$ and $m$ are minimal then this relation will be called the \emph{minimal relationship} between $K_1$ and $K_2$.

On can see that $\DCombp(G)$ is specified as soon as one finds its basic elements and the minimal relationships between each pair of them.
We describe now the algorithm given in \cite{Plumbing} for calculating the basic elements.

Let $K$ satisfy 
\begin{equation}
\label{eq:PartBox}
m(v)+2 \leq \langle K, v\rangle \leq -m(v).
\end{equation}

Construct a sequence of vectors  
$K=K_0,K_1,\ldots,K_n$, where $K_{i+1}$ is obtained from $K_i$ by choosing any vertex $v_{i+1}$ with
$$\langle K_i,v_{i+1}\rangle = -m(v_{i+1}),$$
and then letting
$$K_{i+1}=K_i+2\PD[v_{i+1}].$$
Note that any two vectors in this sequence are equivalent.

This sequence can terminate in one of two ways:
either
\begin{itemize}
\item the final vector $L=K_n$ satisfies the inequality,
\begin{equation}
\label{eq:OtherPartBox}
m(v) \leq \langle L, v\rangle \leq -m(v)-2
\end{equation}
at each vertex $v$ or
\item there is some vertex $v$ for which
\begin{equation}
\label{eq:OtherTermination}
\langle K_{n},v \rangle > -m(v).
\end{equation}
\end{itemize}

It turns out that the equivalence classes in $\DCombp(G)$ which
have no representative of the form $U^m\otimes K'$ with $m>0$ are in
one-to-one correspondence with initial vectors $K$ satisfying
inequality~\eqref{eq:PartBox} for which the algorithm above terminates
in a characteristic vector $L$ satisfying
inequality~\eqref{eq:OtherPartBox}.

\section{The calculation}
 Let us start by computing the basic vectors for the Brieskorn spheres $\Sigma(2, 2n+1, 4n+3)$.
 The plumbing graph G in this case is depicted in the Figure 1. We will write the 
elements of $\Combp(G)$ as row vectors with the 
first four coordinates corresponding to the vertices with weights $-1$, $-2$, $-3$ and $-4n-3$ respectively, and all remaining entries corresponding to $-2$'s on the middle strand ordered by the distance from the root starting with the closest one. 

\begin{figure}[ht!]
\cl{\epsfxsize 1.7in\epsfbox{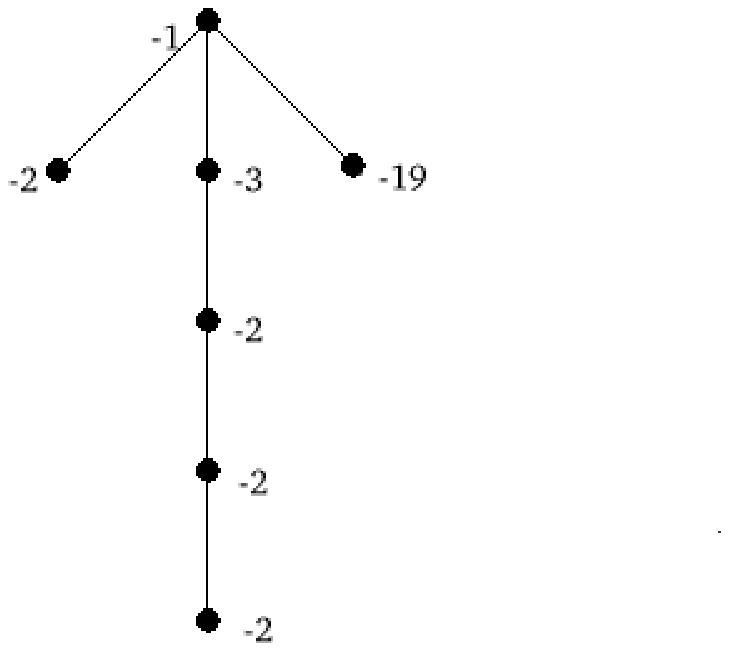}}
\caption{\label{fig:Brieskorn}
{\bf{Plumbing description of $\Sigma(2,9,19)$}}\qua
Here $n=4$; in general the number of $-2$'s on the middle strand is $n-1$.}
\end{figure}

\begin{lemma}
\label{lemma:bv}
For the Brieskorn sphere $\Sigma(2, 2n+1, 4n+3)$ there are $2n$ non-equivalent basic vectors $K_1, K_2, ..., K_{2n}$ where  
$$K_i = (1, 0, -1, -4n-3+2i, 0, 0, ..., 0).$$
\end{lemma}
\begin{proof}
Let $K=(a,b,c,d, ...)$ be a basic vector satisfying~\eqref{eq:PartBox}. Obviously $a=1$. $b$ could be $0$ or $2$, but in the latter case algorithm would stop on the second step with final vector satisfying~\eqref{eq:OtherTermination}. 
For $c$ we again have two possibilities $-1$ and $1$ and the case $c=1$ is eliminated in the same manner on the 4th step.    
The elements corresponding to $-2$'s of the middle strand could be $0$ or $2$. However, if any of them is equal to $2$ then running the algorithm and choosing $v_i$'s from the strand as long as possible we will arrive to a vector of the form $(1, 0 , 1, ...)$ meaning that our initial vector was not basic. To summarize, all the basic vectors are of the form 
$$K=(1,0,-1, d, 0,0,..., 0).$$
 Although the range of $d$ is given by $-4n-1 \leq d \leq 4n+3$, we claim that $K$ is basic only if $-4n-1 \leq d \leq -3$. 

Let us analyze the flow of the algorithm starting with $K=(1,0,-1, d, 0,0,..., 0)$ until the vertex with $m(v) = -4n-3$ is used for the first time or the algorithm terminates. Note that the order in which we choose $v_i$'s does not matter for the outcome of the algorithm, so we decide to extend the sequence as long as possible without using the weight -4n-3 vertex.
This will not give us the whole sequence of the algorithm. We will use some other considerations to figure out if $K$ is a basic vector.
Note that the sequence in which
the vertices are used in this truncated algorithm will depend only on $n$ not on $d$.

 The algorithm starts with the vector $K_0=(1, 0, -1, d, 0, 0, ...0)$ and on the  5th step we get  $K_4=(1,0,-3, d+4,0,0,...0)$. The $v_i$'s used are
 the weight -1, -2 (on the left branch), -1 and -3 vertices, in the given order.

Any consecutive five vectors (among the vectors in the algorithm's sequence) with the first one equal to $(1, 0, -1, d', e, ...)$ and the fifth one being $(1,0,-3, d'+4, e+2, ...)$  using -1, -2, -1 and -3 weighed vertices 
in the order will be called a \emph{cycle}. Note that the entries denoted by "$...$" do not change within a cycle.  

Let \emph{finalizer} to be the sequence of four vectors starting in $(1, 0, -3, d', 0,...,0)$  and
 ending in $(-1, 0, 1, d'+4, 0, ..., 0)$. One uses the weight -1, -2 (on the left branch), -1 vertices in the given order.
Any part of the algorithm's sequence where only the vertices 
 of the long strand with weights -2 are used will be called \emph{transitor}.

We claim that for any positive $n$ the truncated algorithm's sequence consists of $n$ cycles, $n-1$ transitors and the finalizer. Each transitor starts
 with the last vector of the previous cycle and ends with the first vector of the next cycle. The last vertex of the middle strand is used only once in the whole process. The finalizer's first vector is the 5th vector of the last cycle. 

The proof is by induction on $n$. The case of $n=1, 2$ is easily checked by hand. For $n+1$ run the algorithm and choose the same order of the vertices as used for the case of $n$ except for the last ones corresponding to the finalizer. This gives us $n$ cycles and $n-1$ transitors. 
The last cycle will end with the vector $(1,0,-3,x,0,0,...0,2)$ because of the induction hypothesis about the last vertex of the long strand. Now we get a new transitor by choosing $v_i$'s to be the vertices of the long strand in decreasing order of distance from the root. This transitor ends with the vector $(1, 0, -1, x, -2, 0,0,...,0)$. The same vector is the beginning of a new cycle. The end vector $(1,0,-3, x', 0,0,...0)$ of the cycle is the beginning of the finalizer. The finalizer obviously stops the sequence from being continued unless the vertex with weight of -4n-3 is used and this completes the inductive step.

Looking at the way the fourth component of the vectors change one deduces that the truncated algorithm when run for $K=(1,0,-1, d,0,0,...,0)$ stops with $L=(1,0,-1,d+4n+4,0,0,...,0)$. There are three possible cases:
\begin{itemize}
\item one has $d+4n+4<4n+3$. In this case the stopping point of the truncated algorithm coincides with that of the original algorithm. The final vector satisfies the inequalities~\eqref{eq:OtherPartBox} which means that the initial vector is basic. These are the $2n$ vectors given in the lemma. For future reference, we call the sequence of vertices used in this case as the alpha-sequence.    
\item one has $d+4n+4>4n+3$. There is no need to continue the truncated algorithm, because a vector with an entry bigger than the absolute value of the corresponding weight cannot be basic and all vectors in the algorithms sequence are equivalent.
\item case of $d+4n+4=4n+3$. Neither in this case have we a basic vector, because when the truncated sequence is continued after five steps we get a vector with the first entry equal to $3$ which is bigger than the absolute value of the corresponding weight.
\end{itemize}
It follows that all the basic vectors are the ones given in the lemma.
\end{proof}
Now we need to consider the relationships between the basic vectors. 
\begin{lemma}
\label{lemma:rel}
There exist the following relationships between $K_1, K_2, ..., K_{2n}$:

$U \otimes K_1 \sim U^{n+1} \otimes K_2,$

$U \otimes K_2 \sim U^{n} \otimes K_3,$

$U^2 \otimes K_3 \sim U^{n} \otimes K_4,$

$U^2 \otimes K_4 \sim U^{n-1} \otimes K_5,$

\vdots

$U^{n+1} \otimes K_{2n-1} \sim U \otimes K_{2n}.$
\end{lemma}
\begin{proof}
Note that in the relations \ref{equation:rel1} and \ref{equation:rel2} $PD[v]$ is added. Obviously, it is possible to express the same relations with $PD[v]$ subtracted. For each relationship of the lemma there is a sequence of vertices so that subtracting their $PD$'s in order yields relationship along the way.
For brevity, we will refer to the vertices by the column number of the corresponding entry in our vector notation. For example to get the first relationship for each $n$ one starts with $K_1$ and uses the vertex sequence(where we use ";" to make the pattern more visible)
$$ 1,2,3,4,5, ..., n+3,1; 1,2,3,5,...,n+2,1; 1,2,3,5,..., n+1,1; ...;1,2,3,1$$
and follows it by the reversed alpha-sequence. Along the way the relationship $U\otimes K_1 \sim U^{n+1} \otimes K_2$ is obtained. 
To get the second relationship one uses the same sequence but without the first segment and follows the second one by vertex $4$. All other sequences are similar. Note that the alpha-sequence itself ends with $1$.
One can easily see that the first occurrence of "1" and any occurrence of "1,1" in the sequence corresponds to a change in the $U$-power. Given this, it is straightforward if tedious to verify that we indeed get the claimed relationships.
\end{proof}

\begin{lemma}
\label{lemma:length}
The renormalized length $(K_n^2+|G|)/4$ of $K_n$ is equal to zero.
\end{lemma}
\begin{proof}
Let $M$ be the incidence matrix of our graph, with the diagonal elements equal to the weights of the corresponding vertices. Note that 
$$K^2 = K M^{-1} K^T.$$
One verifies by induction that 
$$M^{-1} K_n^T =  (2n, n,n,1,n-1,n-2,...,1)^T.$$
The lemma follows.
\end{proof}

\begin{remark}
\label{remark:l}
Note that the lengths of all basic vectors can be calculated based on the fact that multiplication by $U$ increases the renormalized length by 2
and also one can obtain relationships between any two of the basic vectors using the relationships of the lemma. Among these derived relationships we are interested in those which follow in the most economical way. 
\end{remark}
\begin{lemma}
\label{lemma:min}
The relationships given in the lemma~\ref{lemma:rel} and those which follow from them as in the previous remark are minimal.  
\end{lemma}
\begin{proof}
If the relationships under consideration were not minimal then
$$\Rk \HFred(-Y) < n(n+1)/2.$$ However, since $\HFred$ is supported only in even degrees it follows that $$\Rk \HFred(-Y) = \chi (\HFred(-Y)).$$ By Theorem 1.3 of \cite{AbsGraded} and lemma \ref{lemma:length} this in turn is equal to Casson's invariant of $-Y$, which is seen to be $n(n+1)/2$. 
\end{proof}

\section {Proof of Theorem~\ref{theorem:main}}
We have all the ingredients needed to finish the calculation: the basic vectors and the minimal relationships between them. Note that 
$K_n$ and $K_{n+1}$ have the smallest renormalized lengths. Renormalized length of $K_{i}$  is equal to $(n-i)(n-i+1) = q_{n-i}$ for $i=1,...,n$. One has the minimal relationship 
$$U^{p_i} \otimes K_i \sim U^{p_i+q_{n-i}} \otimes K_{n+1}$$
where $i = 1,2,..., n$.
One just replaces $K_i$ by $K_{n+i}$ and $K_{n+1}$ by $K_n$ in the statements above to get the corresponding results for the remaining vectors. The theorem follows.

\end{document}